\documentclass[11pt]{article}
\input{epsf}
\usepackage{amssymb, latexsym, amsfonts,amsmath, psfrag,
graphics, eepic,colordvi}
\topskip  
\numberwithin{equation}{section}
\newtheorem{theorem}{Theorem}[section]

\newenvironment{kst}
{\setlength{\leftmargini}{2.3\parindent}
\begin{itemize}
\setlength{\itemsep}{-1.1mm}} {\end{itemize}}

\begin{document}
\parskip 6pt

\pagenumbering{arabic}

\def\sof{\hfill\rule{2mm}{2mm}}
\def\ls{\leq}
\def\gs{\geq}

\def\qq{{\bold q}}
\def\MM{\mathcal M}
\def\TT{\mathcal T}
\def\EE{\mathcal E}
\def\OO{\mathcal O}
\def\DD{\mathcal D}
\def\VV{\mathcal V}
\def\NN{\mathcal N}
\def\SS{\mathcal S}
\def\PP{\mathcal P}
\def\lsp{\mbox{lsp}}
\def\rsp{\mbox{rsp}}
\def\pf{\noindent {\it Proof.} }
\def\mp{\mbox{pyramid}}
\def\mb{\mbox{block}}
\def\mc{\mbox{cross}}
\def\qed{\hfill \rule{4pt}{7pt}}
\def\block{\hfill \rule{5pt}{5pt}}
\begin{center}

{\Large\bf     Noncrossing Trees and  Noncrossing
 Graphs}

\vskip 6mm

William Y. C. Chen$^1$  and   Sherry H. F. Yan$^2$

\vskip 3mm

Center for Combinatorics, LPMC, Nankai University, 300071 Tianjin,
P.R. China

{\tt $^1$chen@nankai.edu.cn, $^2$huifangyan@eyou.com}
\end{center}


\vskip 6mm

\noindent {\bf Abstract.} We give a parity
reversing  involution on  noncrossing trees that leads to
 a combinatorial interpretation of a formula  on
noncrossing trees and symmetric ternary trees in answer to
a problem proposed by
Hough.  We use the representation of Panholzer and Prodinger
for noncrossing trees and find
a correspondence between a class of noncrossing trees, called proper noncrossing trees,
and the set of symmetric ternary trees. The second result of this paper
is a parity reversing involution on  connected
noncrossing  graphs which leads to a relation between
 the number of noncrossing trees with
a given number of edges and descents  and  the number of connected
noncrossing graphs with a given  number of vertices and edges.
\medskip




\section{Introduction}

A {\em noncrossing graph} with $n$ vertices is a graph drawn on $n$
points numbered in counterclockwise order on a circle such that
the edges lie entirely within the circle and do not cross each other.
   Noncrossing trees have been
studied by  Deutsch,
Feretic and Noy \cite{DFN},  Deutsch and Noy \cite{DN}, Flajolet and Noy \cite{FN},
Noy \cite{N},  Panholzer and
Prodinger \cite{Pou}. It is well known that the number of  noncrossing trees
with $n$ edges equals the generalized Catalan number $c_n={1\over
{2n+1}}{3n\choose n}$.

\begin{figure}[h,t]
\begin{center}\label{s}
\begin{picture}(50,15)
\setlength{\unitlength}{3mm} \linethickness{0.4pt}
\put(1,1){\circle{2}}\put(0.5,0.5){\small$T$}\put(1,2){\line(2,1){4}}\put(5,4){\circle*{0.2}}
\put(5,4){\line(0,-1){2}}\put(5,1){\circle{2}}\put(4.5,0.5){\small$R$}\put(5,4){\line(2,-1){4}}
\put(9,1){\circle{2}}\put(8.5,0.5){\small$T'$}
\end{picture}
\end{center}
\caption{$T'$ is the reflection of $T$ and $R$ is symmetric.
}\label{pphi}
\end{figure}
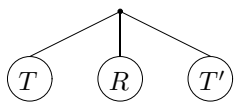

In this paper we are concerned with rooted noncrossing trees.
We assume that $1$ is always the root.
 A {\em descent} is an edge $(i, j)$  such that $i>j$
 and $i$ is on the path from the root $1$ to the vertex $j$.
A {\em ternary tree} is either a single node,
called the root, or it is a root associated with three
ternary trees. A {\em symmetric ternary tree} is a ternary tree
which can be decomposed into a ternary left subtree, a central
symmetric ternary tree and a ternary right subtree that is a
reflection of the left subtree, as shown in Figure \ref{s}.

 Let $\SS_n$ be the set of
symmetric ternary trees with $n$ internal vertices.   A noncrossing
tree is called {\em even} if the number of descents is even.
Otherwise, it is called {\em odd}.  Denote by $\EE_n$ and $\OO_n$
the sets of  even and odd noncrossing trees with $n$ edges,
respectively.
 Let $s_n, e_n, o_n$ be the
cardinalities of the sets $\SS_n, \EE_n, \OO_n$, respectively.
Deutsch, Feretic and Noy \cite{DFN} have shown that
\begin{equation}\label{eq.7}
s_n=\left\{ \begin{array}{ll} {\displaystyle {1\over 2m+1}{3m\choose m}} & \
\mbox{if} \ n=2m,\\
&\\
 {\displaystyle {1\over 2m+1}{3m+1\choose m+1}} & \ \mbox{if} \ n=2m+1.
\end{array}\right.
\end{equation}

 Recently,
Hough \cite{Hou} obtained the generating function for the number
of noncrossing trees with $n$  edges  and a prescribed number of descents.
He also derived the following relation:
 \begin{equation}\label{eq.1}
e_n-o_n=s_n.
\end{equation}

Hough \cite{Hou}  asked the natural question of finding a
combinatorial interpretation of  the above identity (\ref{eq.1}).
In this paper, we obtain a parity reversing involution on
noncrossing trees that leads to a combinatorial interpretation of
(\ref{eq.1}).

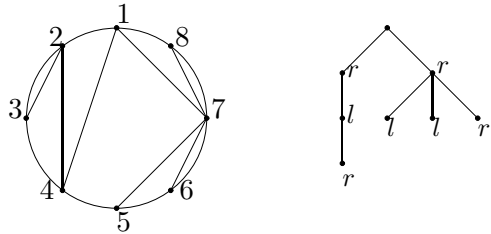
\begin{figure}[h,t]
\begin{center}
\begin{picture}(50,30)
\setlength{\unitlength}{6mm} \linethickness{0.4pt}
\put(2,2){\circle{4}}\put(2,4){\circle*{0.1}}\put(2,4.1){$1$}
\put(2,0){\circle*{0.1}}\put(2,-0.5){$5$}\put(0,2){\circle*{0.1}}
\put(-0.4,2){$3$}\put(4,2){\circle*{0.1}}\put(4.1,2){$7$}\put(0.8,3.6){\circle*{0.1}}\put(0.5,
3.6){$2$}\put(3.2,3.6){\circle*{0.1}}\put(3.3,
3.6){$8$}\put(0.8,0.4){\circle*{0.1}}\put(0.3,
0.2){$4$}\put(3.2,0.4){\circle*{0.1}}\put(3.4, 0.2){$6$}
\put(2,4){\line(-1,-3){1.2}}\put(0.8,0.4){\line(0,1){3.2}}\put(0.8,3.6){\line(-1,-2){0.8}}
\put(2,4){\line(1,-1){2}}\put(4,2){\line(-1,-2){0.8}}\put(4,2){\line(-1,-1){2}}\put(4,2){\line(-1,2){0.8}}

\put(7,1){\circle*{0.1}}\put(7,0.5){\small$r$}\put(7,1){\line(0,1){1}}\put(7,2){\circle*{0.1}}\put(7.1,1.9){\small$l$}
\put(7,2){\line(0,1){1}}\put(7,3){\circle*{0.1}}\put(7.1,
2.9){\small$r$}\put(7,3){\line(1,1){1}}\put(8,4){\circle*{0.1}}
\put(8,4){\line(1,-1){1}}\put(9,3){\circle*{0.1}}\put(9.1,
3){\small$r$}\put(9,3){\line(-1,-1){1}}\put(8,2){\circle*{0.1}}\put(8,1.6){\small$l$}
\put(9,3){\line(0,-1){1}}\put(9,2){\circle*{0.1}}\put(9,1.6){\small$l$}
\put(9,3){\line(1,-1){1}}\put(10,2){\circle*{0.1}}\put(10,1.6){\small$r$}
\end{picture}
\end{center}
\caption{Two representations of a noncrossing
tree}\label{fig.1}
\end{figure}

Our combinatorial interpretation of (\ref{eq.1}) relies on the representation of
noncrossing trees introduced by
 Panholzer and Prodinger \cite{Pou}.
Given a noncrossing tree $T$, we may represent it by a plane tree
with each vertex labeled by $L$ or $R$ with the additional
requirement that the root is not labeled, and the children of the
root are labeled by $R$.  Such a $(L, R)$-labeled tree
representation of $T$ is obtained from $T$ (as a rooted tree) by
the following rule: Given any non-root vertex $j$ of $T$, suppose
that $i$ is the parent of $j$. If $i>j$ then the label of the
vertex corresponding to $j$ is labeled by $L$; otherwise, it is
labeled by $R$. These two equivalent representations of
noncrossing trees are illustrated by Figure \ref{fig.1}.   It is
obvious that a descent in the noncrossing tree in the first
representation corresponds to a $L$-labeled vertex in the second
representation.

The second result of this paper is an expression of the number of
noncrossing trees with $n$ edges and $k$ descents in terms of
the number of connected noncrossing graphs with $n+1$ vertices and $k$ edges.
Noncrossing graphs have been extensively studied by  Flajolet and
Noy \cite{FN}. They derived the following formula for the number of connected
noncrossing  graphs with $n+1$ vertices and $k$ edges,  that is,
\begin{equation}
\label{N-nk}
N_{n,k}={1\over n}{3n \choose n+1+k}{k-1\choose n-1}.
\end{equation}

Hough \cite{Hou} found a combinatorial interpretation of the relation
between the  descent generating function of noncrossing trees and
the generating function for connected noncrossing graphs. By
using the  Lagrange inversion
formula he  obtained
the following  explicit formula for the number of
noncrossing trees with $n$ edges and $k$ descents,
\begin{equation}\label{d-nk}
d_{n,k}={1\over n}{n-1+k\choose n-1}{2n-k\choose n+1}.
\end{equation}

As the second result of this paper, we present a parity reversing involution on
  connected  noncrossing  graphs and obtain
  an expression for the
number $d_{n,k}$ in terms of the number $N_{n,m}$.

\section{An involution on noncrossing trees}

In this section, we give a parity reversing involution on noncrossing trees
which leads to a combinatorial interpretation of the relation (\ref{eq.1}).
We use the  representation of noncrossing trees introduced by
Panholzer and Prodinger \cite{Pou}.
Let $T$ be an even noncrossing tree with $n$ edges and $v$ be a
non-root internal node of $T$.
A vertex $v$ is called a {\em proper}
 vertex if it has   an even number of left children but has no
right child.   If
$T$ is odd, that is, $T$ has an odd number of descents,
 then $v$ is said to be {\em proper} if  $v$ has
 an even number of right children but has no left child.
Otherwise, $v$ is said to be {\em improper}. A noncrossing tree is
said to be {\em proper } if every non-root vertex is proper.
Otherwise, it is said to be {\em improper}. It is obvious that
each odd noncrossing tree is improper.  Let us use $\TT_n$  to denote the set
of proper noncrossing trees with $n$ internal nodes and  let $t_n$
denote the cardinality of $\TT_n$.

Let us recall that a plane tree is said to be an {\em even tree}
if each vertex has an
even number of children. Chen \cite{Ch} gives a bijection $\psi$
between  the set of even plane trees with $2n$ edges and the set
of ternary trees with $n$ internal nodes. A similar bijection is  obtained by
Deutsch, Feretic and Noy \cite{DFN}.  Here we give a brief
description of this bijection.
 Suppose that $T$ is an even plane
tree with $2n$ edges. We  use the following procedure to construct
a ternary tree with $n$ internal vertices.
\begin{kst}
\item [Step 1.] Construct two plane trees $T_1$ and $T_2$ based on
$T$:  $T_1$ is the subtree containing the root and the first two
subtrees of $T$, whereas $T_2$ is the subtree of $T$ by removing
the first two subtrees of the root.

\item [Step 2.]  Combine $T_2$ with $T_1$ by joining $T_2$ as the
last subtree of the root of $T_1$.

\item [Step 3.] Repeat the above procedure for all the nontrivial
subtrees (with at least two vertices) of the root.
\end{kst}

Since each   non-root vertex of a proper even noncrossing tree has
only an even number of left children and has no right child, so we
can discard the labels of its children and  represent a proper
tree as a plane tree such that each subtree of the root is an even
tree.  We define a map $\sigma\colon \TT_n\rightarrow \SS_n$ as
follows.

 \noindent{\bf The map $\sigma$:} Let $T$
be a proper even noncrossing  tree. Let $T_1$ be  the first
subtree of the root.
 The map is
defined by a recursive procedure.
\begin{itemize}
\item Step 1. Assign a vertex as the root and let $\psi(T_1)$ be
the first subtree of the root and  its reflection be the third
subtree of the root.

\item Step 2. Let $T_2$ be the subtree  obtained from  $T$ by
deleting $T_1$,
 and  let $\sigma(T_2)$ be the second subtree of the root.
\end{itemize}

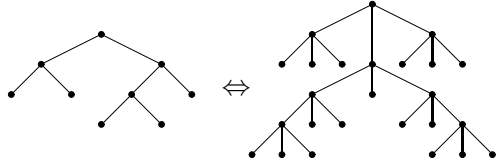
\begin{figure}[h,t]
\begin{center}
\begin{picture}(70,20)
\setlength{\unitlength}{4mm} \linethickness{0.4pt}
\put(1,1){\circle*{0.2}}\put(1,1){\line(1,1){1}}\put(2,2){\circle*{0.2}}
\put(2,2){\line(1,-1){1}}\put(3,1){\circle*{0.2}}\put(2,2){\line(2,1){2}}
\put(4,3){\circle*{0.2}}\put(4,3){\line(2,-1){2}}
\put(6,2){\circle*{0.2}}\put(6,2){\line(-1,-1){1}}\put(5,1){\circle*{0.2}}\put(6,2){\line(1,-1){1}}
\put(7,1){\circle*{0.2}}\put(5,1){\line(-1,-1){1}}\put(4,0){\circle*{0.2}}\put(5,1){\line(1,-1){1}}
\put(6,0){\circle*{0.2}}\put(8,1){$\Leftrightarrow$}
\put(10,2){\circle*{0.2}}\put(10,2){\line(1,1){1}}\put(11,3){\circle*{0.2}}\put(11,3){\line(0,-1){1}}
\put(11,2){\circle*{0.2}}\put(11,3){\line(1,-1){1}}\put(12,2){\circle*{0.2}}\put(11,3){\line(2,1){2}}
\put(13,4){\circle*{0.2}}\put(13,4){\line(0,-1){2}}\put(13,2){\circle*{0.2}}\put(13,4){\line(2,-1){2}}
\put(15,3){\circle*{0.2}}\put(15,3){\line(-1,-1){1}}\put(14,2){\circle*{0.2}}\put(15,3){\line(0,-1){1}}
\put(15,2){\circle*{0.2}}\put(15,3){\line(1,-1){1}}\put(16,2){\circle*{0.2}}
\put(13,2){\line(-2,-1){2}}\put(11,1){\circle*{0.2}}\put(13,2){\line(0,-1){1}}\put(13,1){\circle*{0.2}}
\put(13,2){\line(2,-1){2}}\put(15,1){\circle*{0.2}}\put(11,1){\line(-1,-1){1}}\put(10,0){\circle*{0.2}}
\put(11,1){\line(0,-1){1}}\put(11,0){\circle*{0.2}}\put(11,1){\line(1,-1){1}}\put(12,0){\circle*{0.2}}
\put(10,0){\line(-1,-1){1}}\put(9,-1){\circle*{0.2}}\put(10,0){\line(0,-1){1}}\put(10,-1){\circle*{0.2}}
\put(10,0){\line(1,-1){1}}\put(11,-1){\circle*{0.2}}
\put(15,1){\line(-1,-1){1}}\put(14,0){\circle*{0.2}}\put(15,1){\line(0,-1){1}}\put(15,0){\circle*{0.2}}
\put(15,1){\line(1,-1){1}}\put(16,0){\circle*{0.2}}\put(16,0){\line(-1,-1){1}}\put(15,-1){\circle*{0.2}}
\put(16,0){\line(0,-1){1}}\put(16,-1){\circle*{0.2}}\put(16,0){\line(1,-1){1}}\put(17,-1){\circle*{0.2}}

\end{picture}
\end{center}
\caption{ The map $\sigma$ }\label{fig.3}
\end{figure}

The above  map $\sigma$ is clearly a bijection between $\TT_n$ and $\SS_n$.
    Figure
  \ref{fig.3} is an example.
\begin{theorem}\label{th1.1}
The map $\sigma$ is a bijection between the set of proper noncrossing trees
with $n$  edges and the set
 of symmetric ternary trees with $n$ internal vertices.
\end{theorem}

By using even plane trees as an intermediate structure, we may
obtain a combinatorial interpretation of (\ref{eq.1}) by
constructing an involution on improper noncrossing trees which
 changes the parity of the number of descents.

\begin{theorem} \label{th1.2}
There is a  parity reversing  involution on the set  of  improper
noncrossing trees with $n$ edges. So we have the following relation
\begin{equation}\label{eq.2}
e_n-o_n=t_n.
\end{equation}
\end{theorem}

\pf
 Let $T$ be an improper  noncrossing tree
with $n$ edges. Traverse  $T$ in  preorder and  let
 $v$ be the first encountered  improper node. Define the map
$\phi$ as follows: Case (1), if $T\in \OO_n$ and $v$ has at least
one left child, then $\phi(T)$ is obtained by changing its
rightmost left child
 to a right child and changing all the  children of
 the non-root vertices traversed before $v$ to left children;
Case (2), if $T\in \OO_n$ and $v$ has no left children but has  an
odd number of right children, then $\phi(T)$ is obtained by
changing all the children of $v$ to left children and changing all the
 children of
 non-root nodes traversed before $v$ to left children.

  If $T\in \EE_n$ and $v$ has at
least one right child, then one can reverse the construction in
Case (1). If $T\in \EE_n$ and  $v$ has no right child and has an
odd number of left children, then the construction in Case (2) is
also reversible.
 Hence the map $\phi$
is an involution  on the set of improper noncrossing trees with $n$ edges.
Moreover, one sees that this involution changes the parity of the number of
descents. Thus, we obtain the relation  \eqref{eq.2}. \qed

\begin{figure}[h,t]
\begin{center}
\begin{picture}(50,20)
\setlength{\unitlength}{5mm} \linethickness{0.4pt}

\put(0,1){\circle*{0.2}}\put(0,1){\line(1,1){1}}\put(1,2){\line(1,-1){1}}
\put(2,1){\circle*{0.2}}\put(1,2){\circle*{0.2}}\put(2,1){\line(-2,-1){2}}
\put(0,0){\circle*{0.2}}\put(2,1){\line(-1,-1){1}}\put(1,0){\circle*{0.2}}
\put(2,1){\line(1,-1){1}}\put(3,0){\circle*{0.2}}\put(2,1){\line(2,-1){2}}
\put(4,0){\circle*{0.2}}\put(1,2){\line(1,1){1}}\put(2,3){\circle*{0.2}}
\put(2,3){\line(1,-1){1}}\put(3,2){\circle*{0.2}}\put(2,3){\line(-2,-1){2}}
\put(0,2){\circle*{0.2}}\put(2,3){\line(2,-1){2}}\put(4,2){\circle*{0.2}}
\put(0,0.2){\small$l$}\put(1,0.2){\small$l$}\put(3,0.15){\small$l$}\put(4,0.2){\small$r$}
\put(0,1.2){\small$r$}\put(2,1.2){\small$r$}\put(1,2.2){\small$r$}\put(0,2.2){\small$r$}
\put(3,2.2){\small$r$}\put(4,2.2){\small$r$}
\put(5,1){$\Leftrightarrow$}

\put(7,1){\circle*{0.2}}\put(7,1){\line(1,1){1}}\put(8,2){\line(1,-1){1}}
\put(9,1){\circle*{0.2}}\put(8,2){\circle*{0.2}}\put(9,1){\line(-2,-1){2}}
\put(7,0){\circle*{0.2}}\put(9,1){\line(-1,-1){1}}\put(8,0){\circle*{0.2}}
\put(9,1){\line(1,-1){1}}\put(10,0){\circle*{0.2}}\put(9,1){\line(2,-1){2}}
\put(11,0){\circle*{0.2}}\put(8,2){\line(1,1){1}}\put(9,3){\circle*{0.2}}
\put(9,3){\line(1,-1){1}}\put(10,2){\circle*{0.2}}\put(9,3){\line(-2,-1){2}}
\put(7,2){\circle*{0.2}}\put(9,3){\line(2,-1){2}}\put(11,2){\circle*{0.2}}
\put(7,0.2){\small$l$}\put(8,0.2){\small$l$}\put(10,0.15){\small$r$}\put(11,0.2){\small$r$}
\put(7,1.2){\small$l$}\put(9,1.2){\small$l$}\put(8,2.2){\small$r$}\put(7,2.2){\small$r$}
\put(10,2.2){\small$r$}\put(11,2.2){\small$r$}
\end{picture}
\end{center}
\caption{The involution $\phi$}\label{fig.2}
\end{figure}
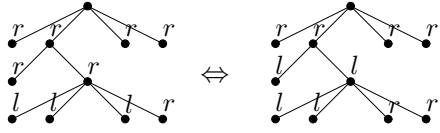

   An  example of the above involution is illustrated in Figure
   \ref{fig.2}.

Combining the bijections in Theorems \ref{th1.1} and \ref{th1.2}, we
get a combinatorial interpretation of the relation
(\ref{eq.1}).  Note that equation (\ref{eq.1}) leads to the
following two combinatorial identities
\begin{eqnarray*}
\sum_{k=0}^{2m-1}(-1)^k{2m-1+k\choose k}{4m-k\choose
2m+1} & =& {2m \over {2m+1}}{3m\choose m}, \\[10pt]
\sum_{k=0}^{2m}(-1)^k{2m+k\choose k}{4m+2-k\choose
2m+2}& = & {3m+1\choose m+1}.
\end{eqnarray*}

\section{An involution on connected noncrossing  graphs}

In this section, we aim to establish a connection between the
number of noncrossing trees with $n$ edges and $k$ descents and the
number of connected noncrossing graphs with $n+1$ vertices and $m$ edges.
\begin{theorem}\label{th.2}
We have the following relation
\begin{equation}\label{eq.4}
\sum_{m=n}^{2n-1}(-1)^{m-n}{m-n\choose k}N_{n,m}=(-1)^kd_{n,k}.
\end{equation}
\end{theorem}

\begin{figure}[h,t]
\begin{center}
\begin{picture}(80,30)
\setlength{\unitlength}{6.2mm} \linethickness{0.4pt}
\put(2,2){\circle{4}}\put(2,4){\circle*{0.1}}\put(2,4.1){$1$}
\put(2,0){\circle*{0.1}}\put(2,-0.5){$4$}\put(4,2){\circle*{0.1}}\put(4.1,2){$6$}\put(0.8,3.6){\circle*{0.1}}\put(0.5,
3.6){$2$}\put(3.2,3.6){\circle*{0.1}}\put(3.3,
3.6){$7$}\put(0.8,0.4){\circle*{0.1}}\put(0.3,
0.3){$3$}\put(3.2,0.4){\circle*{0.1}}\put(3.4, 0.2){$5$}
\put(2,4){\line(0,-1){4}}\put(0.8,3.6){\line(0,-1){3.2}}\put(2,4){\line(-1,-3){1.2}}
\put(2,0){\line(-3,1){1.2}}\put(2,0){\line(1,3){1.2}}\put(3.2,
3.6){\line(0,-1){3.2}}\put(3.2,
3.6){\line(1,-2){0.8}}\put(3.2,0.4){\line(1,2){0.8}}
\put(9,2){\circle{4}}\put(9,4){\circle*{0.1}}\put(9,4.1){$1$}
\put(9,0){\circle*{0.1}}\put(9,-0.5){$4$}\put(11,2){\circle*{0.1}}\put(11.1,2){$6$}\put(7.8,3.6){\circle*{0.1}}\put(7.5,
3.6){$2$}\put(10.2,3.6){\circle*{0.1}}\put(10.3,
3.6){$7$}\put(7.8,0.4){\circle*{0.1}}\put(7.3,
0.3){$3$}\put(10.2,0.4){\circle*{0.1}}\put(10.4, 0.2){$5$}
\put(9,4){\line(0,-1){4}}\put(7.8,3.6){\line(0,-1){3.2}}
\put(9,0){\line(-3,1){1.2}}\put(9,0){\line(1,3){1.2}}\put(10.2,
3.6){\line(0,-1){3.2}}\put(10.2, 3.6){\line(1,-2){0.8}}

\end{picture}
\end{center}
\caption{The canonical spanning}\label{fi.5}
\end{figure}
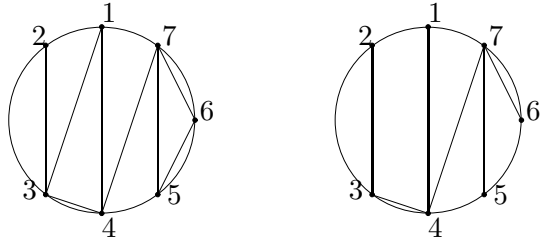

Let $G$ be a connected noncrossing  graph with vertex set $\{1, 2, \ldots, n+1 \}$.
We may construct a unique spanning tree of $G$, which
is called  the {\em canonical spanning tree} of $G$. This
construction can be viewed as a reformation of the traversal procedure  of Hough \cite{Hou}.
Since $G$ is noncrossing,  any cycle of $G$ can be represented by a
sequence $(i_1, i_2, \ldots, i_k)$ such that $i_1 < i_2 < \cdots < i_k$, and
$(i_1, i_2)$, $(i_2, i_3)$, $\ldots$, $(i_{k-1}, i_k)$ and $(i_k, i_1)$ are the
edges of the cycle.
For each cycle $(i_1, i_2, \ldots, i_k)$ in this form, we delete the
edge $(i_1, i_2)$ until we obtain a spanning tree.   An
example is shown in Figure \ref{fi.5}.

Conversely, given a noncrossing tree $T$ with $n$ edges and a
subset $S$ of its descents, we can construct a  connected noncrossing
graph by using the bijection of Hough \cite{Hou} which can be
described  as follows: For each descent $(i,j)$ in $S$,  find the maximal
path of consecutive descents from $j$ back to the root, and let
the first vertex on this path be $v$. From the neighbors of the vertices
on the path from $v$ to $i$ except for the vertices on the path,
choose the neighbor $w$ as the largest vertex less than $j$;
Then add the new edge $(w,j)$ to $T$. We call the new edge $(w,j)$
 the {\em
companion edge} of the descent $(i,j)$.

An edge in $G$  is said to be {\em free} if it is not in the
canonical spanning tree $T$. A descent $(i,j)$ in the canonical
spanning tree of a connected  noncrossing  graph is said to be
{\em saturated } if its companion edge is contained in the
connected  noncrossing graph. Otherwise, it is said to be {\em
unsaturated}.

We now need to consider connected noncrossing graphs in which some
of the free edges are marked. Denote by $\NN_{n,m,k}$ the set of
connected noncrossing  graphs with $n+1$
 vertices and  $m$ edges and $k$ marked free edges.
  It is clear to see that the cardinality of the set
$\NN_{n,m,k}$ is given by
\[ {m-n\choose k}N_{n,m}.\]
Denote by $\NN_{n,k}$ the set of  connected noncrossing graphs
with $n+1$ vertices and $k$ marked free edges.  A descent $(i,j)$
in  the canonical spanning tree of a connected noncrossing
 graph is said to be {\em marked } if its companion edge
 is marked.  Denote by $\DD_{n,k}$ the set of connected noncrossing
graphs  with $n+1$ vertices and  $n+k$ edges such that each descent in
its spanning tree is marked. It follows that
$|\DD_{n,k}|=d_{n,k}$. We will be concerned with the set
$\NN_{n,k}-\DD_{n,k}$, that is, the set of
connected noncrossing graphs with $n+1$ vertices and
$k$ marked free edges  which contain at least one unmarked descent.

Note that two descents $(i, j)$ and $(i',j')$ can not share an end vertex,
namely, $j\not=j'$. A descent $(i,j)$ is said to be smaller than a descent $(i',j')$
if $j<j'$. We now give an involution on the set $\NN_{n,k}-\DD_{n,k}$ that
reverses the parity of the number of free edges.

\begin{theorem}\label{th.1}
There is an involution on  the set
 $\NN_{n,k}-\DD_{n,k}$ that reverses the parity of the number of free edges.
\end{theorem}

\pf Let $G$ be a  connected noncrossing graph in
$\NN_{n,k}-\DD_{n,k}$ with $m-n$ free edges. We define a map
$\psi$ as follows. First, find the minimum unmarked descent
$(i,j)$. We have two cases. Case 1:
  The descent $(i,j)$ is saturated in $G$. We delete the
   companion edge of $(i,j)$ to
get a connected  noncrossing
graph with $n+1$ vertices, $m-n-1$ free edges  and $k$ marked free
 edges.  Case 2: The descent $(i,j)$ is not saturated in $G$.
 We  add the companion edge of $(i,j)$ to get
a connected  noncrossing  graph  with $n+1$ vertices, $m-n+1$ free
edges  and $k$ marked free edges. The operations in the two cases
clearly constitute an involution that changes the number of free
edges by one. \qed

As a consequence of Theorem \ref{th.1}, we obtain  the  identity
(\ref{eq.4}).

To conclude this paper, we remark that
 Theorem \ref{th.2} can be also deduced from the formulas (\ref{N-nk}) and (\ref{d-nk})
for $N_{n,k}$ and $d_{n,k}$ and the following
identity
\begin{equation}
\sum_{m=n}^{2n-1}
(-1)^{m-n-k}{3n \choose n+1+m}{m-1\choose n-1}{m-n\choose k}
={n-1+k\choose n-1}{2n-k\choose n+1}
\end{equation}
that can be verified by using the Vandermonde
convolution \cite[p. 8]{Riordan}
 $${n-m\choose k}=\sum_{i+j=k}(-1)^{i}{m+i-1\choose
i}{n\choose j}.
$$

 \vskip 5mm

\noindent{\bf Acknowledgments.} This work was  supported by
the 973 Project on Mathematical Mechanization,  the
National Science Foundation, the Ministry of Education, and the
Ministry of Science and Technology of China.

%

\end{document}